# 'To be, or not to be, that is the Question': Exploring the 'pseudorandom' generation of texts to write Hamlet from the perspective of the 'Infinite Monkey Theorem'

*Ergon Cugler de Moraes Silva*

Getulio Vargas Foundation (FGV)
University of São Paulo (USP)
São Paulo, São Paulo, Brazil

contato@ergoncugler.com
www.ergoncugler.com

**Abstract**

This article explores the theoretical and computational aspects of the 'Infinite Monkey Theorem', investigating the number of attempts and the time required for a set of pseudorandom characters to assemble and recite Hamlet's iconic phrase, 'To be, or not to be, that is the Question'. Drawing inspiration from Émile Borel's original concept (1913), the study delves into the practical implications of pseudorandomness using Python. Employing Python simulations to generate excerpts from Hamlet, the research navigates historical perspectives and bridges early theoretical foundations with contemporary computational approaches. A set of tests reveals the attempts and time required to generate incremental parts of the target phrase. Utilizing these results, growth factors are calculated, projecting estimated attempts and time for each text part. The findings indicate an astronomical challenge to generate the entire phrase, requiring approximately 2.68×10e69 attempts and 2.95×10e66 seconds — equivalent to 8.18×10e62 hours or 9.32×10e55 years. This temporal scale, exceeding the age of the universe by 6.75×10e45 times, underscores the immense complexity and improbability of random literary creation. The article concludes with reflections on the mathematical intricacies and statistical feasibility within the context of the Infinite Monkey Theorem, emphasizing the theoretical musings surrounding infinite time and the profound limitations inherent in such endeavors. And that only infinity could write Hamlet randomly.

## 1. Introduction

> *"(...) umvvd kqVfv tntlv jWGqf jaeXe **To be** qiaLm abunx oHmzx Zwoua sDadL (...)"*, **or not to be, that is the Question** (Own elaboration, 2024).

In pursuit of answering the intriguing research question regarding the time it would take for a set of pseudorandom characters to recite Hamlet, this article delves into the unique characteristics of pseudorandom values using Python.

By employing Python to generate excerpts from Hamlet, we aim to explore the implications of the Infinite Monkey Theorem. This theorem, originally introduced by Émile Borel in "Mécanique Statistique et Irréversibilité" (Borel, 1913), postulates the idea that given infinite time, a monkey randomly pressing keys on a typewriter would eventually produce a specific text, such as Hamlet (**Attachment 01**). We will shed light on the practical aspects of this concept, examining the role of pseudorandomness in simulating such a literary endeavor.

Building upon the foundational work of Borel, this article also considers subsequent references to the Infinite Monkey Theorem. A. S. Eddington, in 1929, expanded the analogy by proposing that an army of monkeys could potentially write all the books in the British Museum through random typewriting (Eddington, 1929). Sir Jeans, in 1930, extended the discussion, emphasizing the vastness of space and the likelihood of stars forming planetary systems over time, drawing a parallel to the monkey and typewriter scenario (Jeans, 1930).

Arthur Koestler, in 1972, linked the Infinite Monkey Theorem to Neo-Darwinism, describing it as carrying materialism to its extreme limits (Koestler, 1972). Subsequent interpretations by various authors, including Robert M. Pirsig (1974), Douglas Adams (1979), and Scott Adams (1989), introduced humor and further considerations regarding the improbable nature of such random processes.

As we navigate through these historical perspectives, this article aims to contribute to the understanding of the Infinite Monkey Theorem by leveraging modern computational tools. By simulating pseudorandom typing with Python, we will analyze the practical feasibility of this literary experiment, offering insights into the role of chance and pseudorandomness in creative endeavors. The overarching goal is to bridge the theoretical foundations laid by early proponents of the Infinite Monkey Theorem with contemporary computational approaches, exploring the boundaries between randomness and literary creation.

## 2. Materials and methods

In the pursuit of estimating outcomes guided by the principles of the Infinite Monkey Theorem, this paper employs Python methodologies. It involves a simulation where the generation of up to five pseudorandom characters is averaged over ten iterations. In essence, the analysis encompasses the efforts and time invested in producing excerpts such as 'T', 'To', 'To ', 'To b', and 'To be'. These endeavors are extrapolated through a structured loop, iterated a total of ten times. The overarching objective is to project the attained averages as a geometric progression, ultimately leading towards the generation of the complete phrase, 'To be, or not to be, that is the Question'.

### 2.1. Codes

The primary framework for building the dashboard is **Dash**, which is imported from the dash library. Dash enables the creation of dynamic web applications with Python, and it's used here for constructing the layout with HTML components (html). Additionally, the **plotly.express** library is imported as px. The manipulation and handling of data are facilitated by the **pandas** library, imported as pd. The code also incorporates basic Python libraries such as **random** and **string** for generating random characters. The **time** library is used for tracking elapsed time during certain processes. Below, Table 01 presents the initial functionalities.

**Table 01. Tests and code**

| Usage description | Code description |
|---|---|
| This Python code simulates the Infinite Monkey Theorem. The function **generate_partial_text** creates partially random texts of specific sizes. The main function, test_monkey_theorem, takes a target text and iterates over its sizes. It generates random texts until a match with the target is found, recording attempts and elapsed time. The usage section runs the function ten times for the target text 'To be', printing generated texts, attempts, and elapsed time. The code mimics the idea of a monkey typing randomly to have a specific text. | ```python
def generate_partial_text(current_size):
    return ''.join(random.choice(string.ascii_letters + ' ') for _ in range(current_size))

def test_monkey_theorem(target_text):
    generated_text = ''

    for current_size in range(1, len(target_text) + 1):
        start_time = time.time()
        attempts = 0

        while generated_text != target_text[:current_size]:
            attempt = generate_partial_text(current_size)
            elapsed_time = time.time() - start_time
            attempts += 1

            if attempt == target_text[:current_size]:
                generated_text = attempt
                print(f"'{generated_text}' | {attempts + 1:,} attempts | {elapsed_time:.3f} seconds")

# Usage
for i in range(10):
    target_text = "To be"
    attempts = test_monkey_theorem(target_text)
    print("\n")
``` |

**Source:** Own elaboration (2024).

Table 01 results in metrics of attempts and seconds to generate 'T', 'To', 'To ', 'To b' and 'To be'. In this sense, Tables 02 and 03 present these results obtained in the process of analysis, which are used later for projections.

**Table 02. Tests with attempts**

| attempts | 'T' | 'To' | 'To ' | 'To b' | 'To be' |
|---|---|---|---|---|---|
| test 01 | 25 | 10,229 | 38,827 | 8,955,744 | 42,829,442 |
| test 02 | 53 | 4,789 | 113,541 | 4,847,370 | 303,572,721 |
| test 03 | 24 | 5,736 | 441,804 | 6,807,107 | 608,799,906 |
| test 04 | 42 | 715 | 88,579 | 3,679,657 | 127,959,953 |
| test 05 | 89 | 25 | 539,995 | 2,240,406 | 147,686,698 |
| test 06 | 10 | 2,955 | 50,155 | 7,904,165 | 51,964,918 |
| test 07 | 96 | 1,240 | 48,652 | 894,774 | 1,158,891,356 |
| test 08 | 35 | 3,224 | 12,451 | 28,259,478 | 437,532,894 |
| test 09 | 26 | 1,704 | 84,944 | 1,797,972 | 429,123,405 |
| test 10 | 197 | 394 | 172,787 | 15,580,549 | 145,448,105 |
| **average** | **60** | **3,101** | **159,174** | **8,096,722** | **345,380,940** |

**Source:** Own elaboration (2024).

**Table 03. Tests with seconds**

| seconds | 'T' | 'To' | 'To ' | 'To b' | 'To be' |
|---|---|---|---|---|---|
| test 01 | 0.000 s | 0.020 s | 0.081 s | 25.119 s | 136.459 s |
| test 02 | 0.000 s | 0.008 s | 0.240 s | 13.809 s | 965.445 s |
| test 03 | 0.000 s | 0.010 s | 0.977 s | 18.906 s | 1,934.305 s |
| test 04 | 0.000 s | 0.001 s | 0.216 s | 10.124 s | 406.273 s |
| test 05 | 0.000 s | 0.000 s | 1.173 s | 6.451 s | 469.469 s |
| test 06 | 0.000 s | 0.005 s | 0.111 s | 21.793 s | 165.245 s |
| test 07 | 0.000 s | 0.002 s | 0.101 s | 2.341 s | 3,681.907 s |
| test 08 | 0.000 s | 0.006 s | 0.027 s | 77.481 s | 1,389.464 s |
| test 09 | 0.000 s | 0.004 s | 0.285 s | 5.130 s | 1,363.543 s |
| test 10 | 0.000 s | 0.001 s | 0.385 s | 42.391 s | 462.892 s |
| **average** | **0.000 s** | **0.006 s** | **0.360 s** | **22.355 s** | **1,097.500 s** |

**Source:** Own elaboration (2024).

Once the values are generated for the five characters, the averages obtained from ten tests are recorded as variables, as seen in Table 04, below. Furthermore, Table 05 already presents the section of code responsible for running the estimates.

**Table 04. Sets and code**

| Usage description | Code description |
|---|---|
| Here the means estimated by ten previous tests are presented as parameters for geometric progression estimates. | ```python
# Provided Data
attempts_per_part = [60, 3101, 159174, 8096722, 345380940]
times_per_part = [0.0001, 0.0060, 0.3600, 22.3550, 1097.5000]
target_text = "To be, or not to be, that is the Question"
``` |

**Source:** Own elaboration (2024).

**Table 05. Projections and code**

| Usage description | Code description |
|---|---|
| This code calculates growth factors for attempts and time based on the provided **attempts_per_part** and **times_per_part** lists. It then estimates the number of attempts and time for each part of the target text using a geometric progression. The resulting estimates are stored in lists (**text_parts**, **estimated_attempts_f**, **estimated_times**, and **estimated_times_hours**). Finally, a Pandas DataFrame is created to display the estimated data, including text parts, estimated attempts, and estimated time in both seconds and hours. The growth factors are calculated by summing the ratio of each value to its preceding value in the respective lists and dividing by the length of the lists minus one. These growth factors are then used to project the estimated attempts and time for each part of the target text. | ```python
attempts_growth_factor = sum(attempts_per_part[i] / attempts_per_part[i - 1] for i in range(1, len(attempts_per_part))) / (len(attempts_per_part) - 1)
time_growth_factor = sum(times_per_part[i] / times_per_part[i - 1] for i in range(1, len(times_per_part))) / (len(times_per_part) - 1)
text_parts = []
estimated_attempts_f = []
estimated_times = []
estimated_times_hours = []

for i in range(1, len(target_text) + 1):
    text_part = target_text[:i]

    if i <= len(attempts_per_part):
        estimated_attempts = attempts_per_part[i - 1]
        estimated_time = times_per_part[i - 1]
    else:
        estimated_attempts = estimated_attempts * attempts_growth_factor
        estimated_time = estimated_time * time_growth_factor
    text_parts.append(text_part)
    estimated_attempts_f.append(estimated_attempts)
    estimated_times.append(estimated_time)
    estimated_times_hours.append(estimated_time / 3600)

df = pd.DataFrame({
    'Text Part': text_parts,
    'Estimated Attempts': estimated_attempts_f,
    'Estimated Time (seconds)': estimated_times,
    'Estimated Time (hours)': estimated_times_hours})
display(df)
``` |

**Source:** Own elaboration (2024).

Utilizing these metrics, growth factors are calculated, enabling the projection of estimated attempts and time for each text part. The results are showcased in a Pandas DataFrame, providing insights into the intriguing process of simulating the Infinite Monkey Theorem and allowing for further exploration and analysis of the generated data.

## 3. Exploring the results

Table 06, below, presents the results obtained from the geometric progressions, estimating the number of attempts and the estimated time to generate the phrase 'To be, or not to be, that is the Question' from characters in calculations pseudorandoms in Python.

**Table 06. Generated and estimated values**

| Format | Text Part | Estimated Attempts | Est. Time (seconds) | Est. Time (hours) |
|---|---|---|---|---|
| Generated | T | 6,00E+01 | 1,00E-04 | 2,78E-08 |
| | To | 3,10E+03 | 6,00E-03 | 1,67E-06 |
| | To | 1,59E+05 | 3,60E-01 | 1,00E-04 |
| | To b | 8,10E+06 | 2,24E+01 | 6,21E-03 |
| | To be | 3,45E+08 | 1,10E+03 | 3,05E-01 |
| Estimated | To be, | 1,70E+10 | 6,34E+04 | 1,76E+01 |
| | To be, | 8,34E+11 | 3,67E+06 | 1,02E+03 |
| | To be, o | 4,10E+13 | 2,12E+08 | 5,89E+04 |
| | To be, or | 2,01E+15 | 1,22E+10 | 3,40E+06 |
| | To be, or | 9,89E+16 | 7,08E+11 | 1,97E+08 |
| | To be, or n | 4,86E+18 | 4,09E+13 | 1,14E+10 |
| | To be, or no | 2,39E+20 | 2,36E+15 | 6,57E+11 |
| | To be, or not | 1,17E+22 | 1,37E+17 | 3,80E+13 |
| | To be, or not | 5,76E+23 | 7,90E+18 | 2,19E+15 |
| | To be, or not t | 2,83E+25 | 4,57E+20 | 1,27E+17 |
| | To be, or not to | 1,39E+27 | 2,64E+22 | 7,33E+18 |
| | To be, or not to | 6,84E+28 | 1,53E+24 | 4,24E+20 |
| | To be, or not to b | 3,36E+30 | 8,82E+25 | 2,45E+22 |
| | To be, or not to be | 1,65E+32 | 5,10E+27 | 1,42E+24 |
| | To be, or not to be, | 8,11E+33 | 2,94E+29 | 8,18E+25 |
| | To be, or not to be, | 3,99E+35 | 1,70E+31 | 4,73E+27 |
| | To be, or not to be, t | 1,96E+37 | 9,84E+32 | 2,73E+29 |
| | To be, or not to be, th | 9,62E+38 | 5,69E+34 | 1,58E+31 |
| | To be, or not to be, tha | 4,73E+40 | 3,29E+36 | 9,13E+32 |
| | To be, or not to be, that | 2,32E+42 | 1,90E+38 | 5,28E+34 |

| Format | Text Part | Estimated Attempts | Est. Time (seconds) | Est. Time (hours) |
|---|---|---|---|---|
| | To be, or not to be, that | 1,14E+44 | 1,10E+40 | 3,05E+36 |
| | To be, or not to be, that i | 5,61E+45 | 6,35E+41 | 1,76E+38 |
| | To be, or not to be, that is | 2,76E+47 | 3,67E+43 | 1,02E+40 |
| | To be, or not to be, that is | 1,35E+49 | 2,12E+45 | 5,89E+41 |
| | To be, or not to be, that is t | 6,65E+50 | 1,23E+47 | 3,40E+43 |
| | To be, or not to be, that is th | 3,27E+52 | 7,08E+48 | 1,97E+45 |
| | To be, or not to be, that is the | 1,61E+54 | 4,09E+50 | 1,14E+47 |
| | To be, or not to be, that is the | 7,89E+55 | 2,37E+52 | 6,57E+48 |
| | To be, or not to be, that is the Q | 3,88E+57 | 1,37E+54 | 3,80E+50 |
| | To be, or not to be, that is the Qu | 1,90E+59 | 7,90E+55 | 2,20E+52 |
| | To be, or not to be, that is the Que | 9,36E+60 | 4,57E+57 | 1,27E+54 |
| | To be, or not to be, that is the Ques | 4,60E+62 | 2,64E+59 | 7,33E+55 |
| | To be, or not to be, that is the Quest | 2,26E+64 | 1,53E+61 | 4,24E+57 |
| | To be, or not to be, that is the Questi | 1,11E+66 | 8,82E+62 | 2,45E+59 |
| | To be, or not to be, that is the Questio | 5,45E+67 | 5,10E+64 | 1,42E+61 |
| | To be, or not to be, that is the Question | 2,68E+69 | 2,95E+66 | 8,18E+62 |

**Source:** Own elaboration (2024).

**To generate the entire phrase "To be, or not to be, that is the Question", it would be necessary around $2,68 \times 10^{69}$ attempts** (2,680,000,000,000,000,000,000,000,000,000,000,000,000,000,000,000,000,000,000,000,000,000); $2,95 \times 10^{66}$ seconds, $8,18 \times 10^{62}$ hours, or $9.32 \times 10^{55}$ years. It is important to highlight that $9.32 \times 10^{55}$ years is approximately $6.75 \times 10^{45}$ times greater than the estimated age of the universe (which is $1,38 \times 10^{10}$ years).

Therefore, to write "To be, or not to be, that is the Question", it would take $6.75 \times 10^{45}$ times the age of the universe. The following Figures illustrate how the geometric progression almost linearizes the previous values, as it ascends exponentially as the number of values to be generated in the paseudorandomization of characters increases.

**Figure 01. From the beginning to 'To be'**

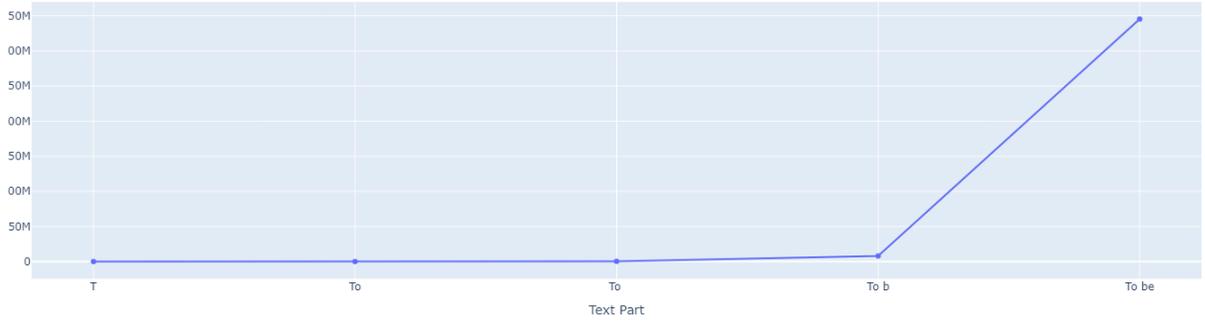

**Source:** Own elaboration (2024).

**Figure 02. From 'To be, or not to be, that is the Quest' to the end**

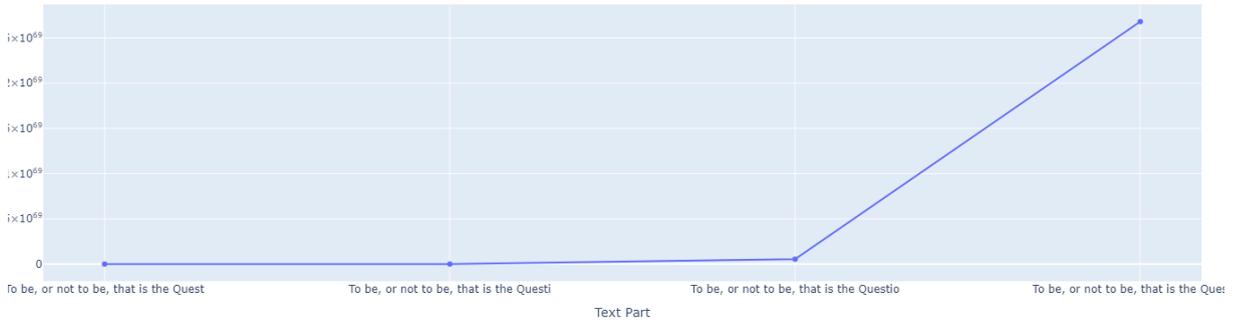

**Source:** Own elaboration (2024).

**Figure 03. The entire data estimation**

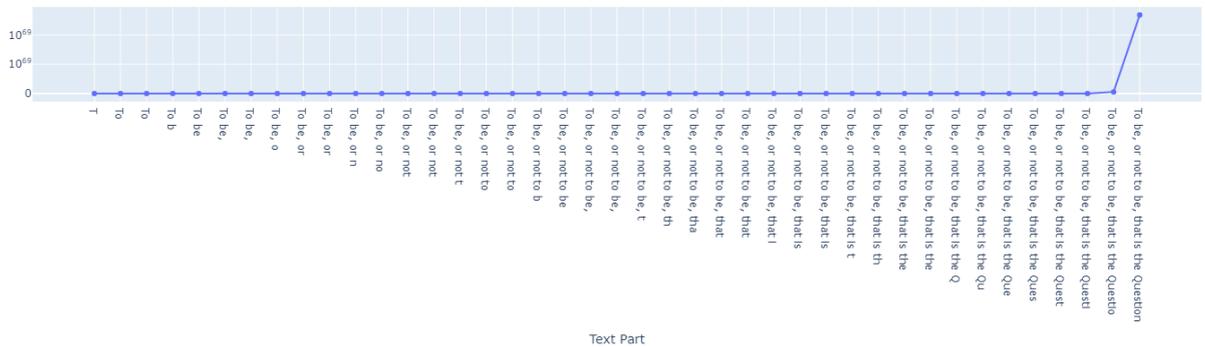

**Source:** Own elaboration (2024).

## 4. Reflections

Drawing from the quantitative insights and the earlier visual representations, it is pertinent to delve into the mathematical intricacies that articulate the likelihood of such

occurrences, fortifying the nuances of this calculated probability. In essence, the probability of randomly assembling characters to form the iconic phrase 'To be, or not to be, that is the Question' is quantified at $4.404^{-71}$, given the intricacies involved in such arrangement:

$$\left(\frac{1}{52}\right)^{41} = 4,404^{-71}$$

If we extrapolate to the 1520 characters that Hamlet (Folio 1, 1623) has (**Attachment 01**), the chance is much smaller, with $4.730^{-2609}$, being:

$$\left(\frac{1}{52}\right)^{1520} = 4,730^{-2609}$$

Just to illustrate, **Attachment 02** shows the representation of values in addition to scientific notation. These calculations highlight the immense improbability of producing specific literary works solely through random character generation, emphasizing the complexity and uniqueness inherent in such texts. However, it also underscores the significance of temporal periods akin to the infinite timeframe posited by the Infinite Monkey Theorem, especially when dealing with extensive text sizes.

This theoretical investigation suggests that given an infinite amount of time, even highly improbable events, such as the random generation of complex literary pieces, become statistically feasible. Thus, while the current probabilities may seem infinitesimally low within the confines of our limited timeframe, the notion of infinite time opens up intriguing possibilities within the context of such theoretical musings.

In the endeavor to generate the entire phrase 'To be, or not to be, that is the Question', the magnitude of the challenge becomes staggering. **It would necessitate approximately $2.68 \times 10^{69}$ attempts**, consuming a colossal time span of $2.95 \times 10^{66}$ seconds, equivalent to $8.18 \times 10^{62}$ hours or $9.32 \times 10^{55}$ years. To put this temporal scale into perspective, it's crucial to emphasize that $9.32 \times 10^{55}$ years is roughly $6.75 \times 10^{45}$ times greater than the estimated age of the universe, which stands at $1.38 \times 10^{10}$ years.

In essence, inscribing 'To be, or not to be, that is the Question' through random character generation would require an astronomical duration, exceeding the age of the universe by a staggering factor. The ensuing figures illustrate how the geometric progression nearly linearizes the prior values, ascending exponentially as the number of values to be generated in the pseudorandomization of characters increases. This underscores the profound complexity of the task at hand, casting light on the limitations and extraordinary timescales involved in such theoretical endeavors. That is, only infinity could write Hamlet randomly.

## 6. Author biography

**Ergon Cugler de Moraes Silva** has a Master's degree in Public Administration and Government (FGV), Postgraduate MBA in Data Science & Analytics (USP) and Bachelor's degree in Public Policy Management (USP). He is associated with the Bureaucracy Studies Center (NEB), collaborates with the Interdisciplinary Observatory of Public Policies (OIPP), with the Study Group on Technology and Innovations in Public Management (GETIP) and with the Monitor of Political Debate in the Digital Environment (Monitor USP). São Paulo, São Paulo, Brazil. Web site: https://ergoncugler.com/.

# 7.     Attachment 01 - Hamlet (Folio 1, 1623)

*To be, or not to be, that is the Question:*
*Whether 'tis Nobler in the mind to suffer*
*The Slings and Arrows of outragious Fortune,*
*Or to take Armes against a Sea of troubles,*
*And by opposing end them: to dye, to sleep*
*No more; and by a sleep, to say we end*
*The Heart-ake, and the thousand Naturall ſhockes*
*That Flesh is heyre too? 'Tis a consummation*
*Deuoutly to be wiſh'd. To dye to sleepe,*
*To sleep, perchance to Dream; I, there's the rub,*
*For in that sleep of death, what dreams may come,*
*When we haue ſhufflel'd off this mortall coile,*
*Muſt giue us pause. There's the respect*
*That makes Calamity of ſo long life:*
*For who would beare the Whips and Scornes of time,*
*The Oppreſſors wrong, the poore mans Contumely,*
*The pangs of diſpriz'd Loue, the Lawes delay,*
*The inſolence of Office, and the Spurnes*
*That patient merit of the unworthy takes,*
*When he himſelfe might his Quietus make*
*With a bare Bodkin? Who would theſe Fardles beare*
*To grunt and ſweat vnder a weary life,*
*But that the dread of ſomething after death,*
*The vndiſcouered Countrey, from whoſe Borne*
*No Traueller returnes, Puzels the will,*
*And makes vs rather beare those illes we haue,*
*Then flye to others that we know not of.*
*Thus Conſcience does make Cowards of vs all,*
*And thus the Natiue hew of Resolution*
*Is ſicklied o're, with the pale caſt of Thought,*
*And enterprizes of great pith and moment,*
*With this regard their Currants turne away,*
*And looſe the name of Action. Soft you now,*
*The faire Ophelia? Nimph, in thy Orizons*
*Be all my ſinnes remembred.*

**William Shakespeares (1623)**

*(a total of 1,520 characters)*

# 8. Attachment 02 - Number Representations

**Table 07. Scientific Notations**

| Scientific Notation | Value |
|---|---|
| $\left(\frac{1}{52}\right)^{41} = 4,404^{-71}$<br><br>"To be, or not to be, that is the Question' | 0.00000000000000000000000000000000000000000000000000000000000000000000**4404** |
| $\left(\frac{1}{52}\right)^{1520} = 4,730^{-2609}$<br><br>the entire Hamlet Folio 1 | 0.0000000000000000000000000000000000000000000000000000000000000000000000000000000000000000000000000000000000000000000000000000000000000000000000000000000000000000000000000000000000000000000000000000000000000000000000000000000000000000000000000000000000000000000000000000000000000000000000000000000000000000000000000000000000000000000000000000000000000000000000000000000000000000000000000000000000000000000000000000000000000000000000000000000000000000000000000000000000000000000000000000000000000000000000000000000000000000000000000000000000000000000000000000000000000000000000000000000000000000000000000000000000000000000000000000000000000000000000000000000000000000000000000000000000000000000000000000000000000000000000000000000000000000000000000000000000000000000000000000000000000000000000000000000000000000000000000000000000000000000000000000000000000000000000000000000000000000000000000000000000000000000000000000000000000000000000000000000000000000000000000000000000000000000000000000000000000000000000000000000000000000000000000000000000000000000000000000000000000000000000000000000000000000000000000000000000000000000000000000000000000000000000000000000000000000000000000000000000000000000000000000000000000000000000000000000000000000000000000000000000000000000000000000000000000000000000000000000000000000000000000000000000000000000000000000000000000000000000000000000000000000000000000000000**4730** |

**Source:** Own elaboration (2024).